\documentclass[leqno]{article} 

\usepackage[blocks]{authblk}

\usepackage{amsfonts,amsmath,amssymb}

\usepackage{mathabx}

\usepackage{mathrsfs}
\usepackage{times}

\usepackage{cases}

\usepackage[latin1]{inputenc}
\usepackage[usenames,dvipsnames]{pstricks}
\newtheorem{theorem}{Theorem}[section]

\newtheorem{lemma}[theorem]{Lemma}
\newtheorem{proposition}[theorem]{Proposition}

\newcommand{\var}{\varphi}

\newcommand{\ai}{\alpha}
\newcommand{\weight}{e^{2s\alpha}}
\newcommand{\T}{\tau}
\newcommand{\norm}[1]{\left\Vert#1\right\Vert}
\newcommand{\abs}[1]{\left\vert#1\right\vert}
\newcommand{\set}[1]{\left\{#1\right\}}

\newcommand{\To}{\longrightarrow}
\newcommand{\p}{\partial}
\newcommand{\R}{\mathbb{R}}

\title{Inverse problems for the diffusion equation with one time observation}

\author[1]{Michel Cristofol}
\author[2]{Alexandre Kawano}

\affil[1]{Aix-Marseille Universite, Institut de Mathematiques de Marseille, CNRS, UMR 7373, Marseille, France; Corresponding author.}
\affil[2]{Universidade de  Sao Paulo, Sao Paulo, Brazil}

\date{\today}

\begin{document}

\maketitle

\begin{abstract}
We examine the inverse problem of recovering the potential term of an n-dimensional heat equation  from a single time solution profile. We show that it is possible to uniquely determine  this potential term via a stability inequality using solution measurements taken at a fixed time, assuming the potential is known within an arbitrary subdomain. The approach is based on the Carleman estimate. Additionally, we discuss the optimality of this type of observation.
\end{abstract}

%

\section{Introduction and main result}  
Consider the diffusion equation, a partial differential equation that describes density fluctuations in a material undergoing diffusion within an open and bounded domain $\Omega$  of $\R^n$  with a sufficiently smooth boundary $\Gamma=\p\Omega$. Given $T>0$, we examine the following boundary value problem for the reaction-diffusion equation with homogeneous Dirichlet boundary conditions:
\begin{numcases}{}
y'(x,t)-\textrm{div}(D(x)\nabla y(x,t)) - p(x) y(x,t) = 0& $(x,t)\in Q$, \nonumber \\
y(x,t)=0 & $(x,t)\in\Sigma$,   \label{1.1} \\
y(x,0)=y_0(x) & $x\in\Omega$,   \nonumber
\end{numcases}
where prime stands for the time derivative, $Q=\Omega \times (0,T)$ and $\Sigma=\Gamma\times(0,T)$. Throughout this paper, $t$ and $x=(x_1,\dots,x_n)$ denote the time variable and the spatial variable respectively, and $y$ denotes the solution, is a scalar function, $p \in L^2(\Omega)$ is the potential term. We will assume that the diffusion matrix $D(x)=(d_{ij})_{1\leq i,j\leq n}$,  where the coefficients $d_{ij}(x)$ are smooth on $\overline{\Omega}$ for $1\leq i,j\leq n$, is positive-definite and satisfy for some positive constant $C>0$
\begin{numcases}{}
 C^{-1}\abs{\xi}^2\leq \sum_{i,j=1}^nd_{ij}\xi_i\xi_j\leq C\abs{\xi}^2, & \label{coer} \\
 \quad d_{ij}(x)=d_{ji}(x), & $x\in\overline{\Omega},\,\,\xi=(\xi_1,\dots,\xi_n)\in\R^n$. \nonumber
\end{numcases}
It is well known that system (\ref{1.1}) possesses an unique solution $y$ such that
\begin{equation}\label{1.2}
y\in H^1(0,T;L^2(\Omega))\cap L^2(0,T;H^2(\Omega)).
\end{equation}

\subsection{ Inverse Coefficient Problem} 

Let $0 < \theta < T$ and let $\omega \Subset \Omega$ be 
an arbitrary fixed subdomain.  Our goal is to determine $p = p(x)$, $x\in \Omega$,
from only the observation data
$$
y(x,\theta), \quad x \in \Omega
$$
assuming that $p(x)$ is known  in $\omega$.
\medskip

We emphasize that the initial condition $y_0(x)$ is unknown. Instead, we assume that measurements of $y(\cdot, \theta)$ can be obtained within in $\Omega$  at any fixed time $\theta>0$. In the two-dimensional case, such observations can be performed, for instance, using thermography at time $\theta$. On the other hand determining initial condition  through direct observation becomes challenging once the heat process has already begun. Therefore, our formulation of the inverse problem holds practical significance.\\
In the field of mathematical and physical inverse problems, the potential identification problem is among the most widely studied. This is because the concept of a potential term naturally corresponds to the reaction term in a linear reaction-diffusion equation, specifically its linear component. Such equations arise in various applications, including hydrology \cite{B72}, heat transfer \cite{BBC85}, population genetics  \cite{AW75}, chemistry \cite{BN91}, biology, and spatial ecology, where they are used to model population dynamics \cite{M02}.

It is well known that the solution $y(x,t)$ and its behavior largely depend on the linear part of the reaction term. In population dynamics, for instance, this term determines whether a population goes extinct, persists, or evolves \cite{CC03}. Consequently, identifying or accurately estimating this coefficient is of critical importance. However, in many physical applications, this term is often unknown or only partially known, as it results from the combined effects of multiple factors and cannot be directly measured. Instead, it is typically inferred through the solution $y(x,t)$, but in most cases, $y(x,t)$ is not available simultaneously for all $x$ and all time $t$.\\
Thanks to Carleman inequalities, the method introduced by Bukgheim and Klibanov \cite{BK81} has enabled the establishment of crucial results, including stability inequalities that link the coefficient to be recovered with the available observations \cite{CR08,IY98,YZ01}. These stability inequalities are particularly important for implementing numerical simulations. The observations considered in these works typically consist of local measurements over $\omega \times (0,T)$, where $\omega$ is a subdomain of $\Omega$ , along with the observation of the solution across the entire domain $\Omega$  at a single time $\theta>0$. On the other hand, many studies on the recovery of the linear part of the reaction term in a reaction-diffusion equation have utilized additional information obtained from a finite number of boundary measurements \cite{DR85,L82, PR86}.\\
Another important inverse parabolic problems are associated to the integral overdetermination. We find  these additive informations in \cite{SY09} associated to parabolic equations parametrized by a diffusion coefficient, and in the monograph \cite{POV00}, in which the authors have established uniqueness results for special classes of coefficients. They work with an equivalent formulation of their inverse problem and they use a Fredholm approach.  In \cite{K05, KP93} the linear parabolic case is addressed and in \cite{BC11}, the authors are interested by the case of Neumann boundary conditions and they study  a certain class of non linear reaction term.\\
The works of Isakov \cite{I91, I06} and Prilepko and Solov'ev \cite{PS87} involve point fixed theorem. In  \cite{I91},  Isakov carried out a stability result requiring additional boundary observations with some sign conditions. In \cite{B95}, Bushuyev considers the case of special non linear source term and in \cite{CY96},  the authors deal with parabolic equations parametrized by a diffusion coefficient and establish a local stability inequality except for a countable set  of parameters.  In \cite{CY08}, the same authors proved a logarithmic stability inequality by the overdetermination at time $T$ assuming analyticity hypothesis of the source term without the knowledge of the initial condition.\\
  In \cite{DYYL09}, the authors applied an optimal control framework to solve an optimization problem in the one-dimensional case and in \cite{H07},  a weak solution approach is used to minimize a coast functional. Additionally, investigations in the case of unbounded domains have been carried out in  \cite{Kawa2} an approach based on transformation into an hyperbolic problem is developed.\\
More recently, a novel approach based solely on pointwise observations of the solution $y(x,t)$ has been applied to obtain uniqueness results in various inverse parabolic problems, including the recovery of potential term in the one-dimensional case \cite{CGHR12,RC10}.\\
 The above list of papers and methods is by no means exhaustive, particularly regarding numerical approaches. However, when summarizing the existing results associated with different types of additional observations, it appears that the most practically relevant supplementary information is the knowledge of $y(x,\theta)$ over $\Omega$. This corresponds to a finite number of direct measurements of the solution.
In other words, selecting this type of overdetermination aligns with a realistic and relatively simple observation method, as it can be obtained using appropriate sensors.\\
Furthermore, we emphasize that the methods relying on Carleman inequalities provide a stability inequality that directly links the coefficient to be reconstructed with a finite number of physical observations. However, stability inequalities derived from Carleman estimates typically require additional observations, such as localized measurements of the solution in both space and time that we avoid in our result. On the other hand, the  method using Carleman inequalities involves the use of weight functions, defined later, associated to a parameter noted $s$ which will play an important role. Indeed, the low or high dependance on $s$ of the  different terms will be crucial to absorb  some of them.\\
In this work, we present a new result that couples a stability inequality with an overdetermination observation at a single time $\theta$ in a multidimensional setting, without requiring knowledge of the initial state. Notably, our result can be formulated in terms of recovering an $x-$dependent term $q(x)$ in a source term $F(x,t) = q(x) k(x,t)$ under the assumption that $k(x,t)$ is known and analytic in time see \cite{BC16}.\\
The final part of this work is devoted to addressing a fundamental question: Is the observation of the solution at a single time $\theta$  sufficient to uniquely determine the potential $p(x)$  in our parabolic problem \eqref{1.1}?

\subsection{Setting and hypothesis}
We define the following spaces:
$$
H^{1,2}(Q)=H^1(0,T;L^2(\Omega))\cap L^2(0,T,H^2(\Omega))
$$
and
$$
\mathcal{C}^{1,0}(\overline{Q})=\mathcal{C}^1([0,T];\mathcal{C}(\overline{\Omega})).
$$
 We fix $C_1>0, C_2>0$ and $C_3>0$  and we further set:\\
  \begin{equation}\label{p}
  \mathcal{A} = \{p \in  H^1(\Omega), C_3 > \norm{p}_{L^2(\Omega)} > C_1\}
 \end{equation} 
   and
      \begin{equation}\label{y_0}
      \mathcal{B} = \{y_0 \in  H_0^1(\Omega) \cap H^2(\Omega), y_0 \geq 0 \mbox{ on } \overline{\Omega} , y_0 \not \equiv 0, \norm{y_0}_{H^2(\Omega)} < C_2\}.
 \end{equation} 
These admissible sets, $\mathcal{A}$ for the potential  and $\mathcal{B}$ for the initial condition come from technical constraints  related to the method used.\\
We consider the two following boundary value problems 
\begin{numcases}{}
y'_1(x,t)-\textrm{div}(D(x)\nabla y_1(x,t)) - p_1(x) y_1(x,t) = 0& $(x,t)\in Q$, \label{1bis} \\
y_1(x,t)=0 & $(x,t)\in\Sigma$, \nonumber
\end{numcases}
and
\begin{numcases}{}
y'_2 (x,t)-\textrm{div}(D(x)\nabla y_2(x,t)) - p_2(x) y_2(x,t) = 0& $(x,t)\in Q$, \label{1ter} \\
y_2(x,t)=0 & $(x,t)\in\Sigma$, \nonumber
\end{numcases}
and we assume that 
\begin{equation}\label{condy2}
\vert y_2(.,\theta) \vert > \delta_0  \mbox{ on } \overline{\Omega} \mbox{ for some } \delta_0>0.
\end{equation}
Note that we can replace this condition by taking the following non homogeneous boundary conditions $y_1(x,t) = y_2(x,t) = h(x,t) $ on $\Sigma$ with $h(x,t) \geq \delta_0$ for some $\delta_0>0$ see \cite{CGR06}. On the other hand, we can remove this condition by adding a suitable control in the system (1.7) as done in \cite{BCGY2009} pp. 697.\\
For simplicity we will rewrite the initial conditions as follows
\begin{equation} 
\left\{
\begin{array}{ll} \nonumber
y_1(x,0) = y_0^1(x), & x \in \Omega,\\
y_2(x,0) = y_0^2(x), & x \in \Omega.
\end{array}
\right.
\end{equation}
By taking the differences between the two previous systems and writing $y= y_1-y_2$ and $p(x) = p_1(x)-p_2(x)$, we get 
\begin{numcases}{}
y' (x,t)-\textrm{div}(D(x)\nabla y(x,t)) - p_1(x) y(x,t) = p(x) y_2(x,t) & $(x,t)\in Q$, \label{sysdif}\\
y(x,t)=0 & $(x,t)\in\Sigma$. \nonumber
\end{numcases}
Let $y$ the solution of (\ref{sysdif}), then by the regularity of the parabolic system (see Evans \cite{Ev} Thm.5, pp.360-361) we have 
$$
y,y'\in H^{1,2}(Q)
$$
moreover there exists $ C>0$ such that
\begin{equation}\label{1.5}
\norm{y}_{L^2(Q)}+\norm{\nabla y}_{L^2(Q)}
+\norm{y'}_{L^2(Q)}\leq C (\norm{y_0}_{H_0^1(\Omega)} + \norm{p}_{L^2(\Omega)})
\end{equation}
and if $y_0 \in H_0^1(\Omega) \cap H^2(\Omega)$
\begin{equation}\label{1.6}
\norm{y(\cdot,t)}_{H^2(\Omega)}\leq C (\norm{y_0}_{H^2(\Omega)} + \norm{p}_{L^2(\Omega)})
\end{equation}
where the constant $C$ depends only on $\Omega$, $T$ and $\norm{D}_\infty$.
\medskip

Further, we consider an arbitrarily fixed function $\hat{p}(x)  \in H^1(\Omega)$.

\subsection{Statement of the main result}

We are now prepared to formally state the main result.
\begin{theorem}\label{T1.1} Let $\omega\Subset\Omega$ be an arbitrary open subdomain of $\Omega$ and $\theta\in (0,T)$. If $y'_1 \in H^{1,2}(Q)$ (resp. $y'_2 \in H^{1,2}(Q)$) satisfies (\ref{1bis}) (resp. (\ref{1ter}),  (\ref{condy2})), then there exists a constant $C = C(\omega,  C_1, C_2, C_3, \theta)> 0$ such that
\begin{equation}\label{1.7}
\Vert p_1-p_2\Vert_{L^2(\Omega)} \le C  \Vert y_1(\cdot, \theta) - y_2(\cdot, \theta)
\Vert_{H^2(\Omega)} 
\end{equation}
for all $(p_1, p_2) \in \mathcal{A}$ satisfying $p_1=p_2= \hat{p}$ in $\omega$.
\end{theorem}

For the corresponding inverse problem in the case of $\theta = 0$, establishing uniqueness and stability using $y_{|\omega_T}$ with an arbitrary $\omega$ remains a highly challenging and open problem. For related inverse problems, see: 
Bukhgeim \cite{Bu}, Isakov \cite{I06}, Lavrent'ev, Romanov and 
Shishat$\cdot$ski\u\i\, \cite{LRS}, Romanov \cite{Rom}.
\medskip

In Section 2, we review the key Carleman estimates from \cite{IY98}. Sections 3 and 4 are dedicated to proving Theorem \ref{T1.1}. Finally, Section 5 explores a conjecture regarding the minimal set of observations of the solution $y(x,t)$  required to determine the potential $p(x)$ without imposing additional assumptions on $p(x)$.
\section{Carleman estimate} 
\setcounter{equation}{0}
Our proof is based on the application of a Carleman estimate, originally introduced in the work of Bukhgeim and Klibanov \cite{BK81}. Additionally, in \cite{K1}, Klibanov established a uniqueness result for a closed problem, albeit without stability, under the assumption of greater regularity for the coefficient to be reconstructed.
\medskip

We set $\omega_T=\omega\times (0,T)$ and we assume that the matrix $D$ satisfies \eqref{coer}. 
Now, we will recall a Carleman estimate for the 
problem:
\begin{numcases}{}
y'(x,t)-\textrm{div}(D(x)\nabla y)(x,t) - p(x) y(x,t)= F(x,t)  & $(x,t)\in Q$, \label{2.1} \\ 
y(x,t)=0\textcolor{red}{,} & $(x,t)\in\Sigma$, \nonumber
\end{numcases}
which is the analogue to the Carleman estimate by Imanuvilov \cite{Im1}.  In a first step, we define appropriate weight functions.
\begin{lemma}[ see \cite{IY98}]\label{L2.1} Let $\omega\Subset\Omega$ be an arbitrary open subdomain, for every open set $\omega_0 \Subset \omega$, there exists a function
$\widetilde{\beta}$ such that: $\widetilde{\beta}\in
\mathcal{C}^2(\overline{\Omega})$, 
$$
\widetilde{\beta} = 0, \quad \p_{\nu}\widetilde{\beta} < 0
\qquad \text{on $\p\Omega$},                      
$$
$$
\vert \nabla\widetilde{\beta} \vert > 0 \qquad
\text{in $\overline{\Omega} \setminus \omega_0$},  
$$
\end{lemma}
We take $K > 0$,
such that 
$$
K \ge 5\max_{\overline{\Omega}} \widetilde{\beta}   
$$
and set
$$
\beta = \widetilde{\beta} + K, \quad
\widehat{\beta} = \frac{5}{4}\max_{\overline{\Omega}}\beta.                                         
$$
Then we introduce the weight functions:
$$
\var(x,t) = \frac{e^{\lambda\beta(x)}}{\ell(t)}, \quad
\ai(x,t) = \frac{e^{\lambda\beta(x)} - e^{\lambda\widehat{\beta}}}
{\ell(t)}, \qquad    \ell(t)=t(T-t)                           
$$
where $\lambda $ is a non negative parameter.\\

Now we recall the first key Carleman estimate.
\begin{lemma}\label{L2.2}
Under the above assumptions, there exists $\lambda_0 > 0$ so that
for any $\lambda > \lambda_0$, there exists a constant 
$s_0(\lambda) > 0$
satisfying the following property: there exists a constant 
$C > 0$, independent of $s$, such that
\begin{multline}\label{2.2}
s^3\int_Q \weight\ell(t)^{-3} \vert y\vert^2 dxdt
+ s\int_Q \weight\ell(t)^{-1} \vert \nabla y\vert^2 dxdt\\
\le Cs^3\int_{\omega_T} \weight \ell(t)^{-3}
\vert y\vert^2 dxdt
+ C\int_Q \weight \vert F \vert^2
dxdt
\end{multline}
for all $y \in H^{1,2}(Q)$ solves (\ref{2.1}) and all $s \ge s_0$.
\end{lemma}
For the application to the inverse problem, we also need to estimate the term $y'$. In particular, to improve the readability of the proof of Lemma \ref{L3.2}, we recall the following Carleman inequality, which can be found in \cite{IY98}.
\begin{lemma}\label{L2.3}
Under the above assumptions, there exists $\lambda_0 > 0$ so that
for any $\lambda > \lambda_0$, there exists a constant 
$s_0(\lambda) > 0$
satisfying the following property: there exists a constant 
$C > 0$, independent of $s$, such that
\begin{equation}\label{2.3}
\frac{1}{s}\int_Q \weight \ell(t) \vert y'\vert^2 dxdt
\le Cs^3\int_{\omega_T} \weight \ell(t)^{-3}
\vert y\vert^2 dxdt 
+ C\int_Q \weight \vert F\vert^2
dxdt
\end{equation}
for all $y \in  H^{1,2}(Q)$ solves (\ref{2.1}) and all $s \ge s_0$.
\end{lemma}
\section{Local estimates for parabolic equation}
\setcounter{equation}{0}
For future reference, we now establish the following technical results, which serve as a type of weighted regularity local estimates for parabolic equations. We consider the following parabolic system:
\begin{equation}\label{3.1}
y'-\textrm{div}(D(x)\nabla y)-p(x) y=F(x,t),\quad x\in\Omega,\,t\in (0,T),
\end{equation}
we assume that $F\in L^2(Q)$ and
\begin{equation}\label{3.2}
F(x,t)=0,\quad x\in\omega,\,\,t\in(0,T).
\end{equation}
\begin{lemma}\label{L3.1}
Let $\gamma \in [0,+\infty)$  and $\omega'\Subset \omega \subset \Omega, \omega'_T = \omega' \times (0,T)$. Then there exists a constant $C=C(\gamma,\omega',\omega)>0$ such that the parabolic estimate
\begin{equation}\label{3.3}
\| e^{ s \alpha} \ell^{-\gamma} \nabla y \|_{L^2(\omega'_T)} \leq C s^{5/4}  \| e^{s \alpha} \ell^{-(\gamma+1)} y \|_{L^2(\omega_T)} ,
\end{equation}
holds for any $y \in H^{1,2}(Q)$ satisfies (\ref{3.1})-(\ref{3.2}) and  for large $s$.
\end{lemma}
Let $\chi  \in \mathcal{C}_0^{\infty}(\Omega)$ be supported in $\omega$ with $\chi(x)=1$ for all $x \in \omega'$ and $\textrm{Supp}(\chi)\subset\omega$.\\
Multiplying $F$ by $e^{2 s \alpha} \ell^{-2\gamma} \chi y$ and integrating over $Q=\Omega \times (0,T)$ we get that
\begin{equation}\label{3.4}
0=\int_{Q}  e^{2s \alpha} \ell^{-2\gamma} \chi(x) F(x,t) y(x,t) dx dt = U+W,
\end{equation}
with
$$
U= \frac{1}{2} \int_{Q}  e^{2s \alpha} \ell^{-2\gamma} \chi(x) \frac{d}{dt}(\abs{y(x,t)}^2) dx dt,
$$
$$
W= -\int_{Q} e^{2s \alpha} \ell^{-2\gamma} \chi(x) \textrm{div}(D(x)\nabla y) y(x,t) dx dt - \int_{Q} e^{2s \alpha} \ell^{-2\gamma} \chi(x) p(x)  y^2(x,t) dx dt.
$$
We handle each of the two terms on the right-hand side of (\ref{3.4}) separately. By integrating by parts in the first term, we obtain:
\begin{equation}\label{3.5}
\left| U\right|=  
\left| \int_{Q} e^{2s \alpha} \ell^{-2 \gamma} \chi(x) (s \alpha_t - \gamma \ell^{-1} (T-2t) )\abs{y}^2 dx dt \right| \leq C s \|  e^{s \alpha} \ell^{-(\gamma+1)} y \|_{L^2(\omega_T)}^2,
\end{equation}
since $| \alpha_t | \leq C \ell^{-2}$. \\
Next we obtain in the same way that
\begin{multline}\label{3.6}
W = \int_{Q} e^{2s \alpha} \ell^{-2 \gamma} \chi \; {}^\T D(x) \nabla y\cdot  \nabla y  dx dt \\
+ \int_{Q} e^{2s \alpha} \ell^{-2 \gamma} ({}^\T\nabla\chi + 2 s \chi {}^\T\nabla\alpha)\cdot D(x)  \nabla y \; y  \, dx dt  \\ 
 -\int_{Q} e^{2s \alpha} \ell^{-2\gamma} \chi(x) p(x) y^2(x,t) dx dt  
=  \sum_{i=1}^4 W_i,
\end{multline}
where ${}^\T v$ stands for the  transpose of $v$.
Further, since
\begin{multline}\label{3.7}
2 \int_{Q} e^{2s \alpha} \ell^{-2 \gamma} \; {}^\T\nabla\chi\cdot D(x) \nabla y \, y \, dx dt
 =\int_{Q} e^{2s \alpha} \ell^{-2 \gamma} \; {}^\T\nabla\chi\cdot D(x) \nabla(\abs{y}^2) dx dt\\
=  -\int_{Q} e^{2s \alpha} \ell^{-2 \gamma} (\textrm{div}(D(x)\nabla\chi)+ 2 s \nabla\chi \cdot\nabla\alpha ) \abs{y}^2 dx dt,
\end{multline}
we get that
\begin{equation}\label{3.8}
\left|  \int_{Q} e^{2s \alpha} \ell^{-2 \gamma} \;  {}^\T\nabla\chi\cdot  D(x)  \nabla y \, y \, dx dt \right| \leq C s \| e^{s \alpha} \ell^{-(\gamma +1 \slash 2)} y\|_{L^2(\omega_T)}^2.
\end{equation}
Moreover, we have
\begin{eqnarray*}
\left| \int_{Q} e^{2s \alpha} \ell^{-2 \gamma} \chi \; {}^\T\nabla\alpha\cdot  D(x) \nabla y \,  y \, dx dt \right| & \leq & s^{-3/2} \| e^{s \alpha} \ell^{-\gamma} \chi^{1 \slash 2} \nabla y \|_{L^2(Q)}^2\\
 &+& s^{3/2}    \| e^{s \alpha} \ell^{-\gamma} \chi^{1 \slash 2} \; {}^\T  \nabla\alpha D(x)  y \|_{L^2(Q)}^2 \\
& \leq & s^{-3/2} \| e^{s \alpha} \ell^{-\gamma} \chi^{1 \slash 2} \nabla y \|_{L^2(Q)}^2 \\
&+& C s^{3/2} \|  e^{s \alpha} \ell^{-(\gamma+1)} y \|_{L^2(\omega_T)}^2.
\end{eqnarray*}
as we have $| \nabla\alpha| \leq C \ell^{-1}$.
Putting this together with (\ref{3.5})-(\ref{3.6})  and \eqref{3.8}, we get:
\begin{eqnarray*}
 W_1 &= &-U-W_2-W_3 - W_4\\
 & \leq & C \left[s \|  e^{s  \alpha} \ell^{-(\gamma+1)} y \|_{L^2(\omega_T)}^2 + s^{-1/2} \| e^{s \alpha} \ell^{-\gamma} \chi^{1 \slash 2} \nabla y \|_{L^2(Q)}^2 \right. \\
 &+& \left. s^{5/2} \|  e^{s \alpha} \ell^{-(\gamma +1)} y \|_{L^2(\omega_T)}^2 \right].
\end{eqnarray*}

Then recalling from (\ref{coer}) that 
$$
W_1 = \int_Qe^{2s \alpha} \ell^{-2\gamma} \chi \; {}^\T D(x) \nabla y  \cdot \nabla y\, dxdt 
 \geq C^{-1} \| e^{s \alpha} \ell^{-\gamma} \chi^{1 \slash 2}   \nabla y  \|^2_{L^2(\omega'_T)},
 $$
and taking large $s$, we end up getting (\ref{3.3}).
We now establish an important  inequality that will enable us to estimate a local weighted norm of $y'$ in terms of a local weighted norm of $y$. The trade-off involves slightly enlarging the domain of integration and carefully managing the weights.
\begin{lemma}\label{L3.2}
Let $\gamma \in [0,+\infty)$  and $\omega'\Subset \omega \subset \Omega$. Then there exists a constant $C=C(\gamma,\omega',\omega)>0$ such that the parabolic estimate
\begin{equation}\label{3.9}
\| e^{s \alpha} \ell^{-\gamma} y' \|_{L^2(\omega'_T)} 
 \leq  C s^2 \| e^{s \alpha}\ell^{-(\gamma+2)} y \|_{L^2(\omega_T)} .
\end{equation}
holds for any $y \in H^{1,2}(Q)$ satisfying (\ref{3.1})-(\ref{3.2}) and for large $s>0$.
\end{lemma}
Choose $\omega''$ such that $\omega' \Subset \omega'' \Subset \omega$. Let $\chi\in\mathcal{C}_0^\infty(\Omega)$ such that $\chi=1$ in $\omega'$ and $\textrm{Supp}(\chi)\subset\omega''$. 
Let $y_*=\chi y$. From the definition of $F$ we have 
$$
y_*'-\textrm{div}(D(x)\nabla y_*) - p(x) y_* =\chi F+\mathcal{Q}_1y \quad \textrm{in}\,\, Q,
$$
where $\mathcal{Q}_1$ is a first order differential operator supported in $\omega''$ and from now for simplicity we will note $\chi F+\mathcal{Q}_1y\equiv \mathcal{Q}_1y \, \textrm{in}\, Q.$\\
 Let $\tilde{y}=\ell^{-(\gamma+1/2)}y_*$. Then $\tilde{y}$ verifies
 $$
 \tilde{y}'-\textrm{div}(D(x)\nabla\tilde{y}) - p(x) \tilde{y} =\ell^{-(\gamma+1/2)}\mathcal{Q}_1y-(\gamma+\frac{1}{2})(T-2t)\ell^{-(\gamma+3/2)}y_*\equiv G
 $$
 Then, applying Lemma \ref{L2.3} with $y$ replaced by $\tilde{y}$, we get that
 \begin{multline}\label{3.10}
\frac{1}{s}\int_Q \weight \ell(t) \vert \tilde{y}'\vert^2 dxdt
\le Cs^3\int_{\omega_T} \weight \ell(t)^{-3}
\vert \tilde{y}\vert^2 dxdt 
+ C\int_Q \weight \vert G\vert^2
dxdt\cr
\leq  Cs^3\int_{\omega_T} \weight \ell(t)^{-2(\gamma+2)}
\vert y\vert^2 dxdt 
+ C\int_Q \weight \vert G\vert^2
dxdt.
\end{multline}
Furthermore
\begin{equation}\label{3.11}
\int_Q \weight \vert G\vert^2 dxdt\leq C\int_{\omega''_T}\weight\ell^{-(2\gamma+1)} \vert \nabla y\vert^2dxdt
+C\int_{\omega''_T}\weight\ell^{-(2\gamma+3)} \vert y\vert^2dxdt.
\end{equation}
Next, applying Lemma \ref{L3.1}, with $\omega'$ replaced by $\omega''$ and $\gamma$ replaced by $\gamma+1/2$, to estimate the first  integral in the RHS.
\begin{equation}\label{3.12}
\int_{\omega''_T}\weight\ell^{-(2\gamma+1)} \vert \nabla y\vert^2dxdt\leq Cs^{5/2}
\int_{\omega_T}\weight\ell^{-(2\gamma+3)} \vert y\vert^2dxdt.
\end{equation}
So, putting (\ref{3.12}) in (\ref{3.11}), we obtain
\begin{equation}\label{3.13}
\int_Q \weight \vert G\vert^2 dxdt\leq Cs^{5/2} \int_{\omega''_T}\weight\ell^{-(2\gamma+3)} \vert y\vert^2dxdt,
\end{equation}
and from \eqref{3.10} and \eqref{3.13} we get 
\begin{equation}\label{3.14bis}
\int_Q \weight \ell(t) \vert \tilde{y}'\vert^2 dxdt \leq C s^4 \int_{\omega_T} \weight \ell(t)^{-2(\gamma+2)}
\vert y\vert^2 dxdt.
\end{equation}
Moreover, since
\begin{equation}\label{3.14}
\int_{\omega'_T}\weight\ell^{-2\gamma} \vert  y'\vert^2dxdt\leq
C \left[\int_{\omega'_T}\weight\ell \vert \tilde{y}'\vert^2dxdt+
\int_{\omega'_T}\weight\ell^{-2(\gamma+1)} \vert y\vert^2dxdt \right]
\end{equation}
we deduce that
\begin{equation}\label{3.15}
\int_{\omega'_T}\weight\ell^{-2\gamma} \vert  y'\vert^2dxdt\leq Cs^4\int_{\omega_T} \weight \ell(t)^{-2(\gamma+2)}
\vert y\vert^2 dxdt.
\end{equation}
The proof is complete.
\section{Proof of the  Theorem \ref{T1.1}} 
\setcounter{equation}{0}

We can now prove Theorem \ref{T1.1} using an argument similar to that of Imanuvilov and Yamamoto \cite{IY98}. In the first step, we establish a new stability inequality, \eqref{4.9}, which involves a weaker norm than the one obtained in \cite{IY98}.\\
For this, we set $Au(x) = \textrm{div}(D(x)\nabla u(x))$ when $u \in 
\mathcal{C}^2(\overline{\Omega})$, and $u\vert_{\p\Omega} = 0$.\\
We recall the system \eqref{sysdif} obtained by  taking the difference between \eqref{1bis} and \eqref{1ter} where $y= y_1-y_2$ and $p(x) = p_1(x)-p_2(x)$. Now we work with  
$$ y'  - A y - p_1 y = p y_2.
$$
Since we can change the scales of $t$, without loss of generality, we
may assume that $\theta = \frac{T}{2}$.
Then the  function $z = y' \in H^{1,2}(Q)$ satisfies
\begin{equation}\label{4.1}
z' = Az + p_1z + p(x) y'_2(x,t) \quad \textrm{in}\,\, Q    
\end{equation}
and
\begin{equation}\label{4.2}
z\left(x, \theta\right) = Ay_\theta + p_1(x) y_\theta +p(x) y_2(x,\theta) \quad \textrm{in}\,\, \Omega,                           
\end{equation}
where we set $y_\theta(x) = y\left(x, \theta\right)$, $x\in \Omega$.
Therefore by Lemmas \ref{L2.2}-\ref{L2.3},  we obtain
\begin{multline}\label{4.3}
s^3\int_Q \weight\ell(t)^{-3} \vert z\vert^2 dxdt+ s\int_Q \weight\ell(t)^{-1} \vert \nabla z\vert^2 dxdt
+ \frac{1}{s}\int_Q \weight \ell(t)\vert z'\vert^2 dxdt\\ 
\le Cs^3\int_{ \omega'_T} \weight \ell(t)^{-3}
\vert z\vert^2 dxdt 
+ C\int_Q \weight\vert p(x) \vert^2 dxdt        
\end{multline}
for all large $s > 0$.\\
By $e^{2s\ai(x,0)} = 0$, 
the Cauchy-Schwarz inequality and (\ref{4.3}), we have
\begin{multline}\label{4.4}
\int_{\Omega} s\left\vert z\left(x, \theta\right)\right\vert^2
e^{2s\ai(x,\theta)}  dx
=  \int^{\theta}_0 \frac{\p}{\p t} 
\left( \int_{\Omega} s\vert z(x,t)\vert^2\weight dx \right) dt \\
= \int_{\Omega}\int_0^\theta 2s^2
(\ai' e^{2s\ai} )z^2dtdx
+ \int_{\Omega}\int_0^\theta
2\left( \frac{1}{\sqrt s} \sqrt{\ell(t)}z'\right)
\left( s\sqrt{s}z\frac{1}{\sqrt{\ell(t)}} \right)\weight dtdx \\
\le  C\int_{\Omega}\int_0^\theta
s^2\ell(t)^{-2} \weight z^2 dtdx
+ C\int_{\Omega}\int_0^\theta
\left( \frac{1}{s}\ell(t) \vert z'\vert^2
+ s^3z^2\ell(t)^{-1} \right) \weight dtdx\\
\le  Cs^3 \int_{\omega'_T}
\weight \ell(t)^{-3} z^2 dxdt
+ C\int_Q \weight \vert p(x)\vert^2 dxdt       
\end{multline}
for all large $s > 0$ and since $ \ell(t)^{-1}< C(T) \ell(t)^{-3}$ .\\
Hence (\ref{4.2}) and (\ref{4.4}) yield
\begin{multline}\label{4.5}
s\int_{\Omega} \vert p(x)\vert^2 \left\vert y_2\left(x,\theta
\right)\right\vert^2 e^{2s\ai(x,\theta)} dx
\le Cs\int_{\Omega} \vert Ay_\theta\vert^2
e^{2s\ai(x,\theta)}dx \\
+ Cs\int_{\Omega} \vert y_\theta\vert^2
e^{2s\ai(x,\theta)}dx\\
+ Cs^3\int_{\omega'_T} \weight \ell(t)^{-3}
\vert z\vert^2 dxdt 
+ C\int_Q \weight \vert p(x) \vert^2 dxdt        
\end{multline}
for all large $s > 0$.
Since $\ai(x,t) \le \ai\left(x, \theta \right)$ for
all $(x,t) \in Q$, we have
\begin{equation}\label{4.6}
\int_Q \weight \vert p(x)\vert^2 dx
\le \int_{\Omega} \vert p(x)\vert^2 \left( 
\int^T_0 e^{2s\ai(x,\theta)} dt \right) dx
= T\int_{\Omega} \vert p(x)\vert^2 
e^{2s\ai(x,\theta)} dx.
\end{equation}
Hence, by the last condition in (\ref{4.6}), 
we can absorb the last term at the right hand side
of (\ref{4.5}) into the left hand side if $s > 0$ is sufficiently large, 
so that
\begin{multline}\label{4.7}
s\int_{\Omega} \vert p(x)\vert^2
e^{2s\ai(x,\theta)} dx
\le Cs\int_{\Omega} \vert Ay_\theta\vert^2 
e^{2s\ai(x,\theta)}dx +  Cs\int_{\Omega} \vert y_\theta\vert^2
e^{2s\ai(x,\theta)}dx \\
+ Cs^3 \int_{\omega'_T} \weight \ell(t)^{-3}
\vert y' \vert^2 dxdt
\end{multline}
for all large $s > 0$.  
By Lemma \ref{L3.2}, we get
\begin{multline}\label{4.8}
s\int_{\Omega} \vert p(x)\vert^2
e^{2s\ai(x,\theta)} dx
\le Cs\int_{\Omega} \vert Ay_\theta\vert^2 
e^{2s\ai(x,\theta)}dx + Cs\int_{\Omega} \vert y_\theta\vert^2
e^{2s\ai(x,\theta)}dx \\
+ Cs^7 \int_{\omega_T} \weight \ell(t)^{-7}
\vert  y \vert^2 dxdt.
\end{multline}
Since
$$
\sup_{x\in\Omega} e^{2s\ai(x,\theta)}, \quad
\sup_{(x,t)\in Q} \ell(t)^{-7} e^{2s\ai(x,t)} < \infty,
$$
fixing $s > 0$ sufficiently large, we obtain the following key result:
\begin{proposition}\label{L4.1} Let $\omega\Subset\Omega$ be an arbitrary open subdomain of $\Omega$ and  $\theta\in (0,T)$. If $y \in H^{1,2}(Q)$ satisfies (\ref{sysdif}), then there exists a constant $C = C(\omega, \theta)> 0$ such that
\begin{equation}\label{4.9}
\Vert p\Vert_{L^2(\Omega)} \le C  (\Vert y(\cdot, \theta)\Vert_{H^2(\Omega)} +\norm{y}_{L^2(\omega_T)})
\end{equation}
for all $p \in L^2(\Omega)$ satisfying $p=0$ in $\omega$.
\end{proposition}
We emphasize that in the previous proposition, unlike classical results, the local observation in  $\omega_T$ is performed using a norm in $L^2(\omega_T)$.
In the second step, we recall the following well-known compactness-uniqueness argument result:
\begin{lemma}\label{L4.2}
Let $X$, $Y$ , $Z$ be three Banach spaces, let $\mathcal{A} : X \to Y$ be a bounded injective linear operator with domain $\mathscr{D}(\mathcal{A})$, and let $\mathcal{K}: X \to Z$ be a compact linear operator. Assume that there exists $C_1>0$ such that
\begin{equation}\label{4.10}
\norm{f}_X \leq C_1 \norm{\mathcal{A}f}_Y +\norm{\mathcal{K}f}_Z,\quad  \forall f \in \mathscr{D}(\mathcal{A}).
\end{equation}
Then there exists $C>0$ such that
\begin{equation}\label{4.11}
\norm{f}_X \leq C\norm{\mathcal{A}f}_Y,\quad \forall f \in \mathscr{D}(\mathcal{A}).
\end{equation}
\end{lemma}
Given $\mathcal{A}$ bounded and injective we reason by contradiction by assuming the opposite to (\ref{4.11}). Then there exists a sequence $(f_n)_n$ in $X$ such that $\norm{f_n}_X = 1$ for all $n$ and $\mathcal{A}f_n \to 0$ in $Y$ as $n$ go to infinity. Since $\mathcal{K}: X \to Z$ is compact, there is  a subsequence, still denoted by $f_n$, such that $(\mathcal{K}f_n)_n$ converges in $Z$. Then this is a Cauchy sequence in $Z$, therefore, by applying (\ref{4.10}) to $f_n-f_m$, we get that $\norm{f_n-f_m}_X \to 0$, as $n, m\to\infty$. As a consequence $(f_n)_n$ is a Cauchy sequence in $X$ so $f_n \to f$ as $n\to\infty$ for some $f \in X$. Since
$\norm{f_n}_X = 1$ for all $n$, then  $\norm{f}_X=1$. Moreover we have $\mathcal{A}f_n \to \mathcal{A} f= 0$ as $n \to \infty$, which is a contradiction to the fact that $\mathcal{A}$ is injective.
\noindent We end the proof of  Theorem \ref{T1.1} by defining
$$
X=\set{ p \in L^2(\Omega),\ p(x)=0\ \textrm{in}\ \omega},\qquad Y=H^2(\Omega),\qquad Z=L^2( \omega_T),
$$
and 
$$
\begin{array}{cccc}
\mathcal{A}  : & X&  \To  & Y  \\
   & p & \longmapsto & \mathcal{A} p=y(\cdot,\theta)  \\
\end{array}
\quad;\qquad\quad
\begin{array}{cccc}
\mathcal{K}  : & X & \To & Z \\
  & p & \longmapsto & \mathcal{K} p =y_{| \omega_T}\\
\end{array}
$$
where $y$ denotes the unique solution to (\ref{1.1}). Then it suffices to prove the two following lemma:
\begin{lemma}\label{L4.3}
The operator $\mathcal{A}$ is bounded and injective.
\end{lemma}
In a first step  the boundedness of $\mathcal{A}$ follows readily from \eqref{p}, \eqref{y_0} and (\ref{1.6}). Then $y$ being solution to the system  (\ref{sysdif}), we deduce from the identities $y(\cdot, \theta) = 0$ and $p = 0$ on $\omega$ that $y'(\cdot,\theta)=0$ on $\omega$. Repeating in the same way we get that the successive derivatives of $y$ with respect to $t$ vanish on $\omega \times\set{\theta}$. Since $y$ is solution to some initial value problem with time independent coefficients, it is time analytic so we have $y= 0$ on $\omega \times(0,T)$. Therefore $p = 0$ on $\Omega$ by (\ref{4.9}), and the proof is complete.
\noindent
We remark that  a similar uniqueness result can be found in \cite{K1} for an unbounded domain via local Carleman estimate.
\begin{lemma}\label{L4.4}
$\mathcal{K}$ is a compact operator from $X$ into $Z$.
\end{lemma}
The operator $\mathcal{K}$ being bounded from $X$ to $H^1(\omega_T)$ as we have 
$$
\norm{y}_{L^2(\omega_T)}+\norm{\nabla y }_{L^2(\omega_T) }+ \norm{y'}_{L^2(\omega_T)}\leq  C \norm{p}_{L^2(\Omega)},
$$
from \eqref{p}, \eqref{y_0}  and (\ref{1.5}) with $C = C(\Omega, T, C_1, C_2, \norm{D}_{\infty})$. Then  the result follows readily from the compactness of the injection $H^1(\omega_T)\hookrightarrow Z=L^2(\omega_T)$.
Finally, putting \eqref{4.9}, Lemmas \ref{L4.2}, \ref{L4.3} and \ref{L4.4} together, we end up getting that 
$$
\norm{p}_{L^2(\Omega)}\leq C\norm{\mathcal{A} p}_Y=C\norm{y(\cdot,\theta)}_{H^2(\Omega)},
$$
which  gives Theorem \ref{T1.1}.

\section{Is it possible to improve the set of observation data to get uniqueness ? }
\setcounter{equation}{0}

A natural question arises: For a parabolic operator of type \eqref{1.1}, does the sole observation of the solution at a fixed time  $T>0$  over the entire open set $\Omega$ allow for the unique determination of the potential $p(x)$ on $\Omega$, without any additional hypotheses or data?\\
Many numerical approaches have been developed, but to our knowledge, no uniqueness result has been proved with this sole information. \\
The particularity of the heat operator, which propagates information at infinite speed, makes this problem very difficult. Indeed, at any time $T>0$, there is a loss of memory in the system.\\
To illustrate this difficulty, one can see in \cite{CR08}, thanks to the figure1  and the figure 2, that there is no apparent link between the solution  at a fixed time $T>0$  and the potential to be recovered. Moreover, Isakov in \cite{I91} Section 3,  addressed this question and constructed a counterexample to prove that there is no uniqueness relationship between the solution at a fixed time $T>0$ and the source term.\\
 We followed this way and we found  a counterexample, showing that there is no injection between the solution data at a fixed positive  time $y(.,T)$ on all the set $\Omega$ and the potential $p(x)$ under the conditions mentioned above.


We reconsider equations \eqref{1bis} and \eqref{1ter}, here reproduced for convenience.
\begin{numcases}{} 
		y'_i(x,t)-\textrm{div}(D(x)\nabla y_i(x,t)) - p_i(x) y_i(x,t) = 0& $(x,t)\in Q$, \nonumber\\
		y_i(x,t)=0 & $(x,t)\in\Sigma$,  \nonumber
\end{numcases}
for $i=1,2$.\\
  Suppose that at $t=T$, $y_1(x,T)=y_2(x,T)=f(x)$, $f\in \mathrm{L}^2(\Omega)$ and that $p_1 \neq p_2$ on $\Omega$.

Let $(\lambda_{i,n})_{n\in\mathbb{N}}$ the eigenvalues, repeated with their respective algebraic multiplicity, of the eigenproblems
\begin{numcases}{}
	(\textrm{div}(D(x)\nabla ) + p_i(x)) S_{i,n}(x) = -\lambda_{i,n} S_{i,n}(x) & $x\in \Omega$, \label{eq:EigenProb}\\
	S_{i,n}(x) = 0, & $x\in\Gamma$,  \nonumber
\end{numcases}
$i=1,2$. By the hypothesis about $D$ and $p_i$ we know that $(S_{i,n})_{n\in\mathbb{N}}$ forms a Hilbertian orthonormal basis of $\mathrm{L}^2(\Omega)$ and that $\lambda_{i,n}^{-1} \rightarrow 0$ as $n\rightarrow +\infty$, $i=1,2$. In fact the existence of such a basis and its properties can be obtained from the following proposition.

\begin{proposition}\label{teo:BaseHilbert} 
    The eigenvectors of problem  \eqref{eq:EigenProb}
     form an orthogonal Hilbertian basis $(S_{i,n})_{n\in\mathbb{N}}$ and the corresponding eigenvalues $(\lambda_{i,n})_{n\in\mathbb{N}}\subset  \mathbb{R}_+$ tend to $+\infty$.  
\end{proposition}
    By the properties of $D(x)$ , the existence of its square root is guaranteed. Let us define the following internal product in $H^1_0(\Omega)$.
    \begin{displaymath}
        \langle u, v \rangle_{H^1_0(\Omega)} = \langle D^{1/2}(x) \nabla u, D^{1/2}(x) \nabla v \rangle_{L^2(\Omega)} + \langle p_i^{1/2}(x) u, p_i^{1/2}(x) v \rangle_{L^2(\Omega)}
    \end{displaymath}
    
    For each fixed $v\in H^1_0(\Omega)$, the linear operator $H^1_0(\Omega) \ni u\mapsto \langle u, v \rangle_{L^2(\Omega}$ is bounded. Therefore by the Riesz Theorem, there is a unique $Tu \in H^1_0(\Omega)$ such that 
    \begin{displaymath}
        \langle Tu, v \rangle_{H^1_0(\Omega)} =  \langle u, v \rangle_{L^2(\Omega)}.
    \end{displaymath}
    
    This defines the operator $T: H^1_0(\Omega) \rightarrow H^1_0(\Omega)$. This linear operator is self-adjoint and compact, by the following reasons.
    
    Since 
    \begin{displaymath}
        \begin{split}
            \langle Tu, v\rangle_{H^1_0(\Omega)} = \langle u, v \rangle_{L^2(\Omega)} = \langle v, u \rangle_{L^2(\Omega)} = \langle Tv, u\rangle_{H^1_0(\Omega)} = \langle u, Tv\rangle_{H^1_0(\Omega)} ,
        \end{split}
    \end{displaymath}
    $T:H^1_0(\Omega) \rightarrow H^1_0(\Omega)$ is self-adjoint. 
    
    To show that it is compact, take a weakly convergent sequence $(u_n)_{n\in\mathbb{N}} \subset H^1_0(\Omega)$, $u_n \rightharpoonup u$. Since the inclusion $H^1_0(\Omega) \hookrightarrow L^2(\Omega)$ is compact, $u_n \rightarrow u$ strongly in $L^2(\Omega)$. Recalling that $\langle Tu, v\rangle_{H^1_0(\Omega)} = \langle u, v \rangle_{L^2(\Omega)}$, we have
    \begin{displaymath}
        \|Tu_k -Tu \|_{H^1_0(\Omega)} = \langle u_k -u, T(u_k -u) \rangle_{L^2(\Omega)} \leq \|T\| \|u_k -u\|^2 \rightarrow 0,
    \end{displaymath}
as $n\rightarrow +\infty$. This shows that $(Tu_k)_{k\in\mathbb{N}}$ converges strongly in $H^1_0(\Omega)$, proving that $T:H^1_0(\Omega) \rightarrow H^1_0(\Omega)$ is compact.

By the Spectral Theorem for compact self-adjoint linear operators, $T$ has an orthonormal Hilbertian basis $(S_n)_{n\in\mathbb{N}}$. The corresponding eigenvalues $(\lambda_n)_{n\in\mathbb{N}}$ tend to zero as $n\rightarrow +\infty$.

The eigenvectors of $T$ are the same of the operator $P\doteq \textrm{div}(D(x)\nabla ) + p(x)$.  If $\lambda_n^{-1}$ is the eigenvalue  of $T$ corresponding to $S_n$, then $-\lambda_n$ is  an eigenvalue of $P$ corresponding to the same $S_n$. In fact, suppose that $T S_n = \lambda_n^{-1} S_n$. Then
\begin{displaymath}
    \langle S_n, v \rangle_{L^2(\Omega)} =  \langle TS_n, v \rangle_{H^1_0(\Omega)} = \lambda_n^{-1}  \langle S_n, v \rangle_{H^1_0(\Omega)} = - \lambda_n^{-1}   \langle P S_n, v \rangle_{L^2(\Omega)}, 
\end{displaymath}
for all $v \in H^1_0(\Omega)$. It means that $P S_n= - \lambda_n S_n$.

On the other hand, suppose that  $PS_n = - \lambda_n S_n$.  Then
\begin{displaymath}
    \langle TS_n, v \rangle_{H^1_0(\Omega)} = \langle S_n, v \rangle_{L^2(\Omega)} = -\lambda_n^{-1} \langle P S_n, v \rangle_{L^2(\Omega)} = \lambda_n^{-1} \langle S_n, v \rangle_{H^1_0(\Omega)}, 
\end{displaymath}
for all $v \in H^1_0(\Omega)$. The conclusion is that $T S_n = \lambda_n^{-1} S_n $.

The properties of $(S_n)_{n\in \mathbb{N}}$ and $(\lambda_n)_{n\in \mathbb{N}}$ follow from the Spectral Theorem for compact self-adjoint linear operators.

Now, we are ready to carry out our counterexample. By assuming that the solutions of  \eqref{1bis} and \eqref{1ter} have the form 
$$y_i(x,t)=\sum_{n\in\mathbb{N}}  g_{i,n}(t) S_{i,n}(x),$$
 we obtain for $i=1,2$, 
\begin{displaymath}
	y_{i}(x,t)= \sum_{n\in\mathbb{N}} c_{i,n} \mathrm{e}^{-\lambda_{i,n}t} S_{i,n}(x),
\end{displaymath}
with $(c_{i,n})_{n\in \mathbb{N}} \in \ell^2$.

Now it is sufficient to choose $(c_{i,n})_{n\in \mathbb{N}}$, $i=1,2$, so that:\\
 $$\sum_{n\in\mathbb{N}} c_{i,n} \mathrm{e}^{-\lambda_{i,n}T} S_{i,n}(x) = f(x), $$
  to prove that the knowledge of solely $y(.,T) $ on all the set $\Omega$ cannot give the potential $p(x)$.
This counter example shows that we need  overdetermining data e.g. assume that $p_1(x)=p_2(x)$ for $x\in \omega$ in our case.\\[5mm]



\end{document}